\documentclass[12pt]{article}
\usepackage[numbers]{natbib}
\usepackage{amssymb}
\begin{document}

\title{\bf Scale free SL(2,R) analysis and the Picard's existence and uniqueness theorem
}

\vspace{3cm}

\author{
Dhurjati Prasad Datta\\  
Department of Mathematics,\\
P.O. North Bengal University, 
Darjeeling, Pin: 734430, India \\
email:dp${_-}$datta@yahoo.com} 
\date{}

\maketitle

\vspace{1cm}

\baselineskip = 22pt

\begin{abstract}
The existence of higher derivative discontinuous solutions to a first order ordinary 
differential equation is shown to reveal a nonlinear SL(2,R) structure of analysis in 
the sense that a real variable $t$ can now accomplish changes not only by linear 
translations $t \rightarrow t + h$ but also by inversions $t \rightarrow 1/t$. We show 
that the real number set has the structure of a positive Lebesgue measure Cantor set. We 
also present an extension of the Picard's theorem in this new light. 
\end{abstract}

\begin{center}
{\bf AMS Subj. Classification}: 26E35, 34F05, 46F30\\
{\bf Key Words}: Scale free, Cantor set, Infinitesimals, nonlinear analysis.
\end{center}

\begin{center}
\em{Int. J. Pure and Applied Math. vol 19, 115-127, (2005)}
\end{center}
 
\newpage

\section{Introduction}
The most basic ingredient of real analysis is of course a variable, living in the real 
number line, which is assumed to undergo changes by linear translations. The framework 
of ordinary analysis is therefore essentially {\em linear}, the relevant concepts of 
limit, differentiability etc being formally defined  to satisfy the simple 
transformation properties of the linear translation group. Notice, for example, that the 
limit $t \rightarrow t_0$ means the linear, continuous decrease ( as an action of the 
translation group) of  the distance $|t-t_0|$ to zero and so on. The application of this 
linear analysis in more complex dynamical problems, for instance, necessitates explicit 
breaking of this linear structure by incorporating nonlinear interactions (terms) in the 
underlying differential (and/or integral) equations. The recent investigations, however, 
reveals the surprising existence of some {\em nonlinear dynamical structures} right at 
the heart of the ordinary analysis \cite{dp1,dp2,dp3, dp4}, thereby extending the 
framework of the linear ordinary analysis to the scale free SL(2,R) analysis. This 
extension is based on the realization that the simplest linear, scale free, differential 
equation, viz.,

\begin{equation}\label{sf}
t{\frac{{\rm d}\tau}{{\rm d}t}}=\tau,\,\tau(1)=1
\end{equation}
    
\noindent {\em admits} a generalized class of nonlinear solutions of the form $\tau_g = 
\tau_s (1 + $O($\eta^{2^n})$), having discontinuities in the $2{^n}$th order derivatives 
at $t =1$ \cite {dp4}. Here, $\tau_s=t$ is the standard solution and $t \rightarrow 
t/t_0 = 1+\eta$ denotes a rescaled variable close to 1. The discontinuity in the $2^n$th 
order derivative stems from the rescaling invariance of eq(\ref{sf}) which allows one to 
inject an infinite sequence of arbitrary, {\em infinitesimal} scale factors into the 
standard solution $\tau_s$. Consequently, the small scale variable $\eta$ is raised to a 
random {\em infinitesimal} variable, residing in a nonstandard real line {\bf R} 
\cite{robin}, rather than in the ordinary real line $R$. Indeed, let $t_\pm = 1 \pm 
\eta,\, \eta>0$. Then obviously $\tau_{s-} =\tau _s(t_{-})= 1-\eta$. Now, it is easy to 
verify that this exact solution can as well be derived starting from an initial 
approximation and then obtaining self similar (multiplicative) correction factors 
recursively, viz., $\tau _s = (1/t_{+})\tau _{1}$ when the correction factor satisfies 
the self similar equation  

\begin{equation}\label{ss}
t_{1-}{\frac{{\rm d}\tau _{1-}}{{\rm d}t_{1-}}}=\tau _{1-}
\end{equation}

\noindent in the smaller scale variable $t_{1 \pm}= 1 \pm \eta^2$ in the neighbourhood 
of 1, and so on, since $t_- = {\frac{1}{(1+\eta)(1+\eta^2)(1+\eta^4)\ldots}}$. The 
freedom of a residual \footnote{The said rescaling is residual in the following sense. 
The trivial rescaling $t \rightarrow t/t_0$ shifts a variable near $t_0$ to one near 1. 
The present rescaling reveals, the so far hidden, freedom of injecting a small scale 
uncertainty in the neighbourhood of a point.} rescaling in eq(\ref{ss}) of the form 
$\alpha _0 t_{1-} = t_{1-}^\prime  = 1- \eta^\prime, \, \eta^\prime = \alpha _0(\eta^2 - 
\epsilon _0/\alpha _0), \, \alpha _0 = 1 +\epsilon _0,\, \epsilon _0$ being an 
arbitrarily small real number drawn at random from the interval $\eta^2 \leq \epsilon _0 
<< \eta $, however, induces a small scale modulations ( fluctuations) in $t_{1-}^\prime$ 
and so on to the higher order iterates $t_{n-}^\prime$ generated self similarly from 
eq(\ref{ss}) over smaller and smaller scales. The new, nontrivial solution \cite{dp4} is 
thus obtained as 

\begin{equation}\label{ns}
\tau _N  =C \prod {\frac{1}{t_{n+}^ {\prime}}}, \, \, C = \prod t_{n+}^{\prime}(0)
\end{equation}

\noindent Clearly, the smaller scales reveals a nonlinear structure in the neighbourhood 
of a real number.  Notice, for example, that $\eta _n^\prime = \alpha _n(\alpha _{n-1}^2 
\eta _{n-1}^{\prime\, 2} -\epsilon _n /\alpha _n),\, \epsilon _n$ representing higher 
order uncertainties and $\alpha _n t_{n+} \neq t_{n+}^\prime$. The second inequality 
tells the nontriviality of the residual rescaling symmetry, which is responsible for the 
higher derivative discontinuity of the new solution $\tau _g$ at $t=1$, provided, of 
course, the rescaling is performed an infinite number of times ( for a proof see 
Appendix). Notice also that $\tau _g$ is 2nd derivative discontinuous if the nontrivial 
rescalings are injected at the 2nd and an infinite sequence of subsequent iterates, 4th 
derivative discontinuous when the freedom of this rescaling is used for the first time 
at the 3rd, instead of 2nd, iterate and so on. In general, the new solution would have 
discontinuous $2^{n-1}$ derivative if the nontrivial rescaling is utilized at the $n$th 
iterates (the standard solution $\tau _s$ is recovered when this rescaling is postponed 
till the infinite number of iterations \cite{dp4} (c.f. Appendix)).   

\par In view of this class of discontinuous solutions to the simplest scale invariant 
linear equation (\ref{sf}), the ordinary real number system $R$ needs to be extended, 
accommodating this small scale nondifferentiable structure, to a nonstandard -like real 
number system {\bf R}. Notice that the generalized solution of eq(\ref{sf}) can indeed 
be written as $\tau _g = t (1 \pm \epsilon \phi(t_1)),\, \phi(t) = \tau _N (t)/t,\,t_1 = 
\epsilon t$ where ${\frac{{\rm d}\phi}{{\rm d}\ln t_{1-}}}= 0$,  $\tau _N $ being the 
new solution (\ref{ns}). Consequently, every ordinary real number, for instance, $t_0$  
gets extended to a set of the form $(t_0)_\epsilon = \{t:\, t_0-\epsilon < t < 
t_0+\epsilon, \, \epsilon \neq 0,\, $O($\epsilon^2) =0, \, t \neq t_0 \}\subset${\bf R} 
where $\epsilon$ denotes the irreducible uncertainty even when the real number is 
evaluated with an  ``infinitely  precise measurement". Clearly, the inclusion $\ldots 
(t_0) _ {\epsilon^2} \subset (t_0)_\epsilon$ is satisfied for higher precision 
evaluation of $t_0$. We note that the generalized solution, as it is obtained in 
\cite{dp4}, is constructed purely as a function defined in $R$ and hence the extended 
set {\bf R} constructed out of the generalized solution coincides exactly with $R$ viz., 
{\bf R}= $R$. We remark that this feature of our extension distinguishes it from the 
work of Robinson \cite{robin}. In the traditional approach of nonstandard analysis, 
infinitesimals are defined as an equivalence class of sequences of real numbers. As a 
result, infinitesimals, being constructed out of real numbers, can be thought of as {\em 
extraneous} to the set of ordinary real number set. Our results, on the other hand, 
suggests that {\em infinitesimals are indeed members of the real number set $R$} and so 
might  instead be defined formally  even in the real number set $R$. This approach is 
being developed recently by Bose \cite{mkb}. In view of this existence of infinitesimals 
in $R$, with random properties, the new solution $\tau _g \in$ {\bf R} can as well be 
interpreted to reveal a new mode of change that is available to a (real) variable $t$, 
viz., $t$ can change near an infinitesimal neighbourhood of $t=1$ by (random) 
inversions: $t_- \rightarrow t_-^{-1}=t_+,\, t_\pm = 1 \pm \eta$, for an infinitesimal 
$\eta$, defined by $\eta \neq 0$ but O($\eta^2) =0$, besides the ordinary translations 
over the ordinary real number line.  Notice that the scale free equation (\ref{sf}) is 
invariant under the transformation $t \rightarrow 1/t,\,\tau \rightarrow 1/\tau$ and so 
$t_- \rightarrow t_-^{-1}=t_+$ indeed stands for the generalized solutions considered 
here. The possibility of (random) inversions breaks the exact determinacy of a real 
variable. Further, the scale invariance tells that  the small scale fluctuating 
behaviour of a real variable ( over $\eta$) is reproduced self similarly over the 
smaller nonlinear scales ($\eta^{2^n}$), though the self similarity is respected in a 
statistical sense. Accordingly, the linear framework of the ordinary analysis gets 
extended to a SL(2,R) analysis. One therefore infers that the real number system $R$, as 
it is ordinarily understood in the context of the linear (ordinary) analysis, enjoys a 
host  of richer ``dynamical" properties which get revealed in the SL(2,R) analysis. We 
have already discussed some of these new features of this analysis in \cite{dp3}. Here, 
we present an analysis of the Picard's theorem \cite{ode} in the light of the new 
generalized solutions of eq(\ref{sf}). To this end we first discuss how the arguments 
used in the simple equation (\ref{sf}) can be extended to more general equations (Sec. 
2). Next we discuss some salient features of the computational model introduced in 
\cite{dp4} and also show  that {\bf R} has the structure of a positive Lebesgue measure 
Cantor set \cite{ott}. The (box counting) dimension of a point in {\bf R } ( hence in 
$R$ ) has the value $\nu$, where $ \nu= {\frac{\sqrt 5 -1}{2}}$ is the golden mean 
(Sec.3). Before stating the Picard's theorem in this new light, we present an extension 
of the ordinary Riemann integration  on {\bf R} (Sec.4). 
 
\section{General case}

\par To prepare for a discussion of the Picard's theorem, let us clarify further the set 
up of the scale free SL(2,R) analysis. To this end, we consider the more general linear 
ODE of the form   

\begin{equation}\label{gen}
{\frac{{\rm d}\ln \tau }{{\rm d}\ln t}}= f(t,\tau ),\,\tau(1)=\tau _0
\end{equation}

\noindent Let $\tilde \tau(t)$ be the corresponding standard solution. This solution is 
exactly determinable, in principle, in the linear analysis. In the extended framework 
this, however, corresponds to the zeroth order solution only. Based on this solution, 
one can, however, generate a more {\em complex} solution in the following way. Following 
the steps outlined in Sec.1, let us write the generalized solution close to $t=1$ (more 
precisely in a lhs neighbourhood of 1) as $\tau _{g-}(t_{g-}) \equiv\tau (t_-) = \tilde 
\tau (1/t_+)\tau _-^{\prime}$.  Then it is easy to verify that $\tau _-^{\prime}$ would 
satisfy an equation of the form 

\begin{equation}\label{gen1}
{\frac{{\rm d}\ln \tau _-^{\prime}}{{\rm d}\ln t_{1-}}}= f^{\prime}(t_{1-},\tau 
_-^{\prime})
\end{equation}

 \noindent on the smaller nonlinear scale $\ln (1 - \eta^2)$ where $f^{\prime}(t_{1-}, 
\tau _-^{\prime}) = {\frac{(t_+ f(t_-,\tau _-)- t_-f(1/t_+, \tilde \tau(1/t_+)))}{t_+ - 
t_-}}$. Clearly, the generalized solution $\tau _g$ would belong to the higher order 
nondifferentiable classes of solutions when the infinite sequence of nontrivially 
rescaled (equivalently, random) variables $t_{n-}^{\prime}$ are incorporated as 
indicated above.

\par Notice that the exact self similarity on the  smaller (nonrandom) scales $\ln (1 - 
\eta^2)$ and etc is obtained only for the simplest equation (\ref{sf}). For, a $t-$ 
dependence either of the form $f(t)$ or $f(t,x)$ breaks this exact self similarity. 
Nevertheless, an approximate self similarity is maintained on scales {\em randomised} by  
the nonlinear $t$-dependence, besides the above mentioned freedom of rescaling,  when 
the random scale is defined, for instance, by $  t_{1-}^{\prime} = t_1^{\mu} \approx 1 - 
\mu\eta^2$ where $\mu = f^{\prime}(t_{1-},\tau _-^{\prime})/f(t_{1-},\tau _-^{\prime}) 
\approx 1$, which remains almost constant over the scale $\eta^2$. Note that, 
randomisation of scales by explicit nonlinearity in the ODE already simulates the 
residual rescaling freedom of eq(\ref{sf}).  The emergence of random behaviour in 
dynamical systems (for instance, one dimensional maps, higher order nonlinear ODEs etc) 
from truncation errors are well known   in deterministic chaos. Here, we point out 
equivalent behaviour even in one dimensional {\em linear} ODEs as an effect of 
infinitesimal nonlinear scales (elements) in real number system. As examples we consider 
two simple equations in the following.

\vspace{.8cm}

Case 1. \underline {Modulated exponential}

\vspace{.8cm}

Small scale modulations in the ordinary exponential function defined by 

\begin{equation}\label{exp}
{\frac{{\rm d}\tau}{{\rm d}t}}=\tau
\end{equation}

\noindent is expected in the present formalism. One verifies that the generalized  
solution near $t=1$ is obtained as 

\begin{equation}\label{exp1}
{\rm exp}_{g}(t_-) = e^{1/t_+ + 1/t_{1+} + \ldots}
\end{equation}

\noindent where $t_{1+} = 1 +(1/t_+ - t_-) = 1 + {\frac{\eta^2}{1-\eta}}$ so that $\mu = 
1/(1-\eta)$, and the equality in (\ref{exp1}) is up to a multiplicative factor. In a 
computational problem (truncated) nonlinear terms in $\mu$ would build up over time to 
introduce small scale {\em random} oscillations in the ordinary exponential. Such 
behaviour would  of course be common to {\em any} deterministic functions of linear 
analysis. 

\vspace{.8cm}

Case 2. \underline{ Quadratic nonlinearity}

\vspace{.8cm}

Let us consider the equation

\begin{equation}\label{qad}
t{\frac{{\rm d}\tau}{{\rm d}t}}=\tau^2,\,\tau(1)=1
\end{equation}
          
\noindent having the exact solution $\tau _0 = 1/(1-\ln t)$. Writing $\tau_g(t_{g-}) 
\equiv \tau _0(t_-)   = \tau _{0}(1/t_+)\tau_{1-}$, as usual, we get 

\begin{equation}\label{qad1}
{\frac{{\rm d}\ln \tau_{1-}}{{\rm d}\ln t_{1-}}}=\mu \tau _{1-}
\end{equation}

\noindent where $\mu = {\frac{t_+ - t_-(\frac{\tau_0(1/t_+)}{\tau_0(t_-)})}{t_+ - t_-}} 
{\frac {1-\ln t_{1-}^{\prime}}{1 - \ln t_{-}}} \approx 1$. We point out that an equation 
of the form (\ref{qad1}) is truly nonlinear in the sense that nonlinear random scale 
$t_{1-}^{\prime} = t_{1-}^{\mu} \approx 1-\mu \eta^2$ is itself determined
by the ``unknown" function $\tau_{1-}$ through $\mu$ and hence does not lead to an exact 
solution in the ordinary sense.
\par To conclude this section, we remark that the infinitesimal fluctuations induced to 
the solutions to eq(\ref{gen}) because of SL(2,R) inversions would have a universal 
feature determined by the functional dependence $f(t,\tau)$ close to $t=1$.

\section{Infinitesimals, computation and Cantor sets}
 
In ref \cite{dp4} we presented a computational model of infinitesimally small real 
numbers. Here we explain the relationship of inversions and infinitesimally small 
numbers in the context of this model. Recall that every ordinary real number $t$ is 
replaced by $t_\epsilon = t(1 \pm \epsilon \phi(t_1))$ so that $t_\epsilon$ essentially 
corresponds to a continuum set $(t)_\epsilon$. Now in a computational problem, a real 
number is only treated as a finite precision (decimal/binary) representation. 
Consequently, the real number set is covered by a countable collection of disjoint open 
intervals of the form $t_\epsilon$, where $\epsilon$ now denotes the finite precision 
with which the real numbers are evaluated (calculated). Notice that  $\epsilon$ can have 
values, in decimal representation, of the form 0.5 (for integers, so that 1, for 
instance, stands actually for the open interval (0.5,1.5) and so on), 0.05, 0.005, 
0.0005,$\ldots$. As a result, in any computation the real number set, although conceived 
as a  connected, continuum set (without any gap), is actually realized as a totally 
disconnected, countable set (equivalent to the set of rationals). As one improves upon 
the accuracy and approaches to the infinite precision the connected continuum set of 
linear analysis is thought to have been recovered only in the limit. In view of the 
scale free extension, one, however, meets with an obstacle in the form of {\em 
irreducible infinitesimal uncertainties}, thereby realizing the continuum structure of 
the real number set but for the total disconnectedness. In fact, the extended real 
number set acquires the structure of a Cantor set with a positive Lebesgue measure 
\cite{ott}. 

{\bf{Lemma:}} The set {\bf R} has the box counting dimension $\sigma,\, \sigma = 1 + \ln 
(1 + \epsilon \phi)/\ln t$. Further, $\sigma$ has the asymptotic limits $\lim _{t 
\rightarrow 0^+}\sigma = \nu$ and $\lim _{t \rightarrow \infty}\sigma = 1 +\nu$ where 
$\nu = {\frac{\sqrt 5 -1}{2}}$.

{\em Proof:} In view of the structure of new solution presented in Sec.1, {\bf R} is an 
extension of $R$ in the following sense. Let $N_{0 \delta}$ be the number of open balls 
(intervals) of radius $\delta$ covering an  open interval of $R$. The box counting 
dimension of $R$ is $ \lim _{\delta \rightarrow 0^+}\frac{\ln N_{0\delta}}{\ln 
\delta^{-1}}$  which, in fact, equals 1. To cover {\bf R} the set $\{N_{0\delta}\}$, 
however, is not sufficient. One, in fact, needs to subdivide each of subintervals of 
size $\delta$  into smaller subdivisions of size $\delta^2$ and this process of 
subdivisions is continued recursively over smaller and smaller scales 
$\delta^4,\,\delta^8 \ldots$. In each step {\bf R} is covered better but for a tiny 
residual part, which could only be measured (detected) looking at the next level of 
finer  scale. Accordingly, these subdivisions over finer scales facilitates one to 
measure the required ``length"  of {\bf R} more and more accurately only in a 
progressive manner.\footnote{The concept of a perfect (or total) covering (and/ or 
infinite precision) is, however, illusory. What is more meaningful is the continual 
progress towards perfection \cite{dp1,dp2}. }    Let, in the first step of subdivisions, 
the fractional increase in the number of subdivisions be $N_{1\delta^2}/N_{0\delta} (= 
\lambda <1)$, where $N_{1\delta}^2$ is the number of subdivisions of scale  $\delta^2$ 
needed to cover the boundary ``points" of the original interval. Notice that the 
interior points in the open interval concerned  are already covered well by the zeroth 
level covering balls  and so finer scale covering balls are irrelevant for these points. 
Continuing this process sequentially over smaller scales (as indicated above) one then 
gets, for instance, as the fractional increase  $N_{2\delta^4}/N_{1\delta^2}$ and so on. 
Because of the (relative) self similarity  of the subdivisions   one concludes  
$N_{2\delta^4}/ N_{1\delta^2} = N_{1\delta^2}/N_{0\delta}$ etc. Consequently, total 
number of subdivisions of different sizes needed to cover the said interval of {\bf R} 
is $\tilde N =\sum N_{n \delta^{2^n}} = N_{0\delta} (1 + \lambda +\lambda^2 + \ldots)$. 
Hence the box counting dimension of $\bf R$, evaluated relative to scale $\delta$, is 
obtained as $\sigma = \lim _{\delta \rightarrow 0^+}\frac{\ln \tilde N}{\ln \delta} = 1 
+ \frac{\lambda}{1 - \lambda}$. 

\par Let $\bf t \in \bf R$. Then the ``length" of the open interval $(\bf 0, \bf t)$ 
would be given by $t^{\sigma}$ where $t \in R$. Again, in view of of the new solution we 
write ${\bf t} = t (1 + \epsilon \phi)$. Comparing the two we deduce the first result. 
The value of $\lambda$ is thus determined by the nontrivial factor in ${\bf t}$. Its 
exact value is, however, unimportant for our purpose. That it is nonzero is actually 
more significant.

\par To prove the remaining part of the lemma, we note that $\lim _{t \rightarrow 
0^+}\sigma = \lim _{t^{-1} \rightarrow \infty} (1 - \epsilon (\frac{t^{-1}}{\ln 
t^{-1}})\tau(t_1) + {\rm higher \, order \,terms}) = \lim _{t^{-1} \rightarrow \infty} 
(1 - \epsilon (\frac{t^{-1}}{\ln t^{-1}})\tau(t_1) + \epsilon^2 (\frac{t^{-2}}{\ln 
t^{-2}})(\tau(t_1))^2 + \ldots )  $ ($\tau$ denotes the new solution of (\ref{ns})). As 
$t^{-1}$ approaches to $\infty$ through higher and higher order scales, there exists a 
sufficiently large interval in $R$, when $\pi(t^{-1}) = \frac{t^{-1}}{\ln t^{-1}}\approx 
t^{-1} \approx \epsilon^{-1}$. Accordingly, in this interval $\lim _{t \rightarrow 
0^+}\sigma =  \lim _{t^{-1} \rightarrow \infty}(1/ (1 + \tau(t_1)) ) = \lim _{t^{-1} 
\rightarrow \infty}(1/ (1 + (1/(1 + \epsilon^2 (\frac{t^{-1}}{\ln t^{-1}}) \tau(t_2)) ) 
$ where $t_2 = \epsilon^2 t$ and so on. Notice that the r.h.s of the first equality 
mimics exactly an application of inversion on the first two terms of the logarithmic 
expansion of $\sigma$. The desired limit now follows from the continued fraction 
expansion of the golden mean $\nu$.  Similar arguments also hold for the other limit 
when one makes use instead  $\phi(t) = t \tau(t_1^{-1}),\, \tau(t_1^{-1}) = 1/\tau(t_1)$ 
\cite{dp2}. $\Box$

\vspace{.7cm}

{\bf Remark 1} To see the origin of inversion, let us reconsider the computational model 
introduced above. Let a computation could only distinguish integer numbers, so that $1 
\equiv (1)_{.5}$, for instance. Consequently, $(1)_\epsilon = 1 +\epsilon \mu,\, \mu = 
\mu _-\cup \mu _+ ,\, \mu _-= (-1,0),\, \mu _+ = (0,1),\, \epsilon = 0.5$. The number 1, 
therefore, corresponds to {\em any} number drawn perfectly at {\em random} from 
(0.5,1.5) (the assumption of randomness is codified in the residual rescaling symmetry 
c.f. Sec. 1). This is realized in the above by letting $\epsilon$ instead to lie in the 
interval $(-0.5,0)\cup(0,0.5)$ and $\mu = \phi\sim$O(1). The irreducible fluctuations as 
encoded in the new solution $\phi$ now tells that a point in (-0.5,0) could fluctuate to 
a point in (0,0.5) (by inversion, viz., by scaling and flipping of sign) and vice versa, 
instead of pure translations, as assumed in ordinary analysis.  Further, the especial 
role  of inversions in a Cantor set is also intuitively clear because a point in a 
Cantor set can change to another only by an inversion (i.e., by a discrete jump) because 
of the gaps (voids) separating the Cantor points \cite{dp3}.  

\vspace{.7cm}

{\bf Remark 2} To visualize the Cantor set like structure of $(0)_\epsilon$, we note 
that $(0)_\epsilon = (-\epsilon, 0_0)\cup (0)_{\epsilon^2} \cup (0_0, \epsilon)$, where 
$0_0$ denotes the zero (0) at level 0, and so on recursively over smaller scales. 
Accordingly, in this extended framework, the number 0, for instance, is represented by 
the set $\cup^\infty _0 \{(-\epsilon^{2^n}, 0_n)\cup (0)_{\epsilon^{2^{n+1}}}\cup (0_n, 
\epsilon^{2^n})\}$, $0_n$ being the $n$th level zero. Because of the natural inclusion 
$0_n \supset 0_{n+1}$, the set $\{0_n \}$ represents a finer and finer realization of 
the originally coarse grained zero (0) as the accuracy of the computation is increased. 
In linear analysis, however, these fine structures are ignored so that the ordinary 
singleton $\{0 \}$ is reproduced. Notice that the usual singletons $\{ t_0 \}$ in the 
reduced (coarse grained) set $R$ are reinserted once the nontrivial SL(2,R) generator is 
frozen.

\vspace{.7cm} 

{\bf Remark 3} The significance of $\pi(t^{-1})$ in the prime number theorem  and other 
number theoretic results will be considered elsewhere.  

\section{The Picard's theorem}

\par In view of the Cantor set like structure of {\bf R} (every point of $R$ is actually 
got replaced by a Cantor set, c.f., remark 2), the definition of the ordinary (Riemann) 
integration needs to be extended. Notice that the Cantor set of $(1)_\epsilon $, for 
instance, is written as $(1)_\epsilon = \cup _0^\infty [(1_n - \epsilon^{2^{n}}, \, 1_n 
+ \epsilon^{2^n})- \{1_n \}]$ (c.f., remark 2) of length $2 \sum _0^\infty 
[\epsilon^{2^{n}}- (\epsilon^{2^{n}} -\epsilon^{2^{n+1}})] = 2 \sum _1^\infty \epsilon 
^{2^n} $ which is essentially zero, being of higher order infinitesimal, relative to the 
first order infinitesimal $\epsilon$.  Notice that the length of voids in $(1)_\epsilon$ 
is $2 \sum (\epsilon^{2^{n}} -\epsilon^{2^{n+1}}))= 2 \epsilon$, which equals to the 
length of $(1-\epsilon, 1 + \epsilon)$ in $R$, so that the Lebesgue measure of this 
Cantor set is effectively zero. 

\par  Notice that our iteration process, leading to the new solution of eq(\ref{sf}), 
not only reveals a {\em multiplicative} structure of the nontrivial neighbourhood of 1, 
but also factors the differential operator in (\ref{sf}) into a countably infinite set 
of self similar operators over nonlinear scales $t_{n}^{\prime}$. Accordingly, the 
corresponding integration measure defined on the said Cantor set, incorporating random 
inversions, is defined by the following replacement 

\begin{equation}\label{int} 
\int _{1-\eta}^1 {\rm d}\ln t_- \rightarrow { \mathcal{E_{\bf R}}}\int _{1- \eta}^1 {\rm 
d}\ln t_- \equiv \sum _1^{\infty}\,\int^{t_{n+}^{\prime{-1}}}_1 {\rm d}\ln 
t_{n+}^{\prime}  
\end{equation}

\noindent where $\eta$ is an infinitesimal variable, so that the new solution of 
eq(\ref{sf}) $\tau _g =C \prod {\frac{1}{t_{n+}^ {\prime}}}, \, \, C = \prod 
t_{n+}^{\prime}(0)$ \cite{dp4} is retrieved by direct integration from eq(\ref{sf}). 
Notice that 1 in the integral of the r.h.s. of eq(\ref{int}) corresponds actually to 
$1_n$ with the definition that $\ln 1_n = 0$ for each $n$. More generally, the extended 
integral for a function $f$ in {\bf R} is defined by 

\begin{equation}\label{int1}
{\mathcal{E_{\bf R}}}\int _{\bf a}^{\bf t} f({\bf t}){\rm d}{\bf t} = \int _a^t f(t){\rm 
d}t + \epsilon \sum \int _1^{t^{\prime{-1}}_n} [f(t^{\prime{-1}}_n) + \tau 
_g\,t{\frac{{\rm d}f}{{\rm d}t}}]\,{\rm d}\ln t^{\prime}_n
\end{equation}  

\noindent where we neglect higher order infinitesimals, and use  rescaling invariance   
$ t \rightarrow t/a$ of the logarithmic differential. We also make use of $\tau _g=t 
\frac{{\rm d}\tau _g}{{\rm d}t}$, ${\rm d}\tau _g = \tau _g \sum {\rm d}\ln 
t^{\prime}_n$ and ${\bf t} = t + \epsilon \tau _g$.  Notice that the second term in the 
bracketed integral arises from the first order Taylor's expansion of $f(\bf t)$ and the 
extended integral can be considered as the anti-derivative of the ODE 

\begin{equation}\label{der}
{\frac{{\rm d}F}{{\rm d}{\bf t}}}= f({\bf t})
\end{equation} 

\noindent in {\bf R} (c.f. eq(\ref{gen})). As an example, the modulated exponential 
(\ref{exp1}) is recovered when eq(\ref{exp}) is integrated following the generalized 
integral (\ref{int1}) with $f(t) = 1$.     

\vspace{.5cm}

{\bf Remark 4} Clearly, the infinite sum in (\ref{int1}) is uniformly convergent in 
every closed interval near $t=1$ for continuous $df/dt$. Let us also point out that the 
infinite set of rescaled nonlinear variables $t^{\prime}_n$ are intrinsically 
randomised, as indicated in relation to eq(\ref{gen}), implicitly by $f(t)$ itself. 
Consequently, the extension of integral as defined in eq(\ref{der}) should be 
interpreted as a generalized integral over a set with small scale (infinitesimal) random 
elements. We make a more rigorous treatment of this generalized integral elsewhere. 

\vspace{.5cm}

\par Proceeding at this heuristic level (i.e., forgetting the randomness so that all the 
variables $t$ and $t^{\prime}_n$ are well defined functions of the real variable 
$\eta,\,\,0<\eta<<1$) we now reexamine the Picard's theorem in this new light. 
        
\par  Let us consider eq(\ref{gen}), viz.,

\begin{equation}\label{pic}
{\frac{{\rm d}\ln \tau}{{\rm d}\ln t}}= f(t,\tau),\,\tau(1)=\tau _0
\end{equation}
 
 \noindent where $f$ is assumed to be a $C^\infty$ function, for simplicity. The 
Picard's theorem then guarantees a unique solution in the neighbourhood of $t=1$. It 
also tells that the corresponding solution is $C^\infty$. The proof requires one to 
convert the IVP (\ref{pic}) to the equivalent integral equation

\begin{equation}\label{inteq}  
\ln \tau = \ln \tau _0 + \int _{1}^{t}f(t,\tau){\rm d}\ln t
\end{equation}
 
 \noindent and then to construct a sequence of approximations $\tau _n(t)$ satisfying 
  
\begin{equation}\label{inteq1}  
\ln \tau _n = \ln \tau _0 + \int _{1}^{t}f(t,\tau _{n-1}){\rm d}\ln t
\end{equation}

\noindent converging uniformly in a closed interval within the said neighbourhood to the 
required solution. One verifies that the proof of the Picard's theorem applies also to 
the generalized integral (\ref{int1}) without major modifications. Consequently, we 
state the modified Picard's theorem in the SL(2,R) analysis as follows

{\bf Theorem} Let $f({\bf t}, \tau)$ be $C^{2^n-1}$ in {\bf R}. Then the ODE 

\begin{equation}\label{pic1}
{\frac{{\rm d} \tau}{{\rm d}\bf t}}= f({\bf t},\tau),\,\tau(1)=\tau _0
\end{equation}
 
\noindent has a unique $C^{2^n-1}$ solution in a suitable neighbourhood of ${\bf t} = 
1$.

\vspace{.7cm}

{\bf Remark 5} The function $f$ in the above theorem inherits the differentiable 
structure of $\bf t$. If instead $f \in C^m,\, m < 2^{n}-1 $ then the solution would 
obviously be only $C^m$. 

\section*{Appendix}

The new solution $\tau _N$ reduces to the standard solution $\tau _s$ if the residual 
rescalings are performed only up to a finite number of iterations.  Let the rescaling be 
terminated  at the $n$th iteration, so that $t^\prime _{(n+1)+}= 1 +  \eta _{n+1},\, 
\eta _{n+1} = \alpha _n^2 \eta^{\prime\,2}_n$ \cite{dp4}. Subsequent iterations (without 
rescalings) would then lead to a factor $1 - \eta _{(n+1)}^2$ in $\tau _N$, viz.,

$$\tau _{N-} =C {\frac{1}{t_{+}}}{\frac{1}{t_{2+}^ {\prime}}}\ldots {\frac{1}{t_{n+}^ 
{\prime}}}(1 - \eta _{(n+1)}^2)$$

\noindent Consequently, the uncertainty (or twist, to put it in a pictorial way ) 
introduced in the solution through each nontrivial rescaling and subsequent inversion 
would {\em unwind} gradually, so to speak, leading to the ordinary (nonrandom) $\tau _s$ 
(multiplied by $n$ irrelevant   scale factors).  Things, however, change drastically if 
rescalings are performed an infinite number of times, instead. Indeed,  the statement 
$n$ tends to  infinity means, of course, that $n$ is becoming larger and larger, 
signifying  an {\em unending} process. As a result, the unwinding can never be initiated 
in the iteration process  for an infinitely large $n$.

\end{document}